\documentclass[12pt]{amsart}
\usepackage{amssymb}
\usepackage{geometry}

\theoremstyle{definition}

\theoremstyle{remark}

\numberwithin{equation}{section}

\input{tcilatex}
\begin{document}
\title{J-Groups of Lens Spaces modulo Odd Prime Powers}
\author{Mehmet K\i rdar}
\address{Department of Mathematics, Faculty of Arts and Science, Nam\i k
Kemal University, Tekirda\u{g}, Turkey}
\email{mkirdar@nku.edu.tr}
\date{12 June, 2013}

\begin{abstract}
We describe the J-Groups of lens spaces modulo odd prime powers via
relations and rational representations.
\end{abstract}

\maketitle







\section{Introduction}

$J$-Groups of lens spaces were studied a lot especially by Japanese
mathematicians after J. F. Adam's solution of the vector fields on spheres
problem by making use of the $J$-groups of real projective spaces. In [1],
I. Dibag described the group structure of $J$-Groups of lens spaces modulo
prime powers. His work involves heavy elementary number theory. Another
articles to look at are [8] for odd prime powers and [2] for powers of 2.
The number theory work there is again quite complicated.

The purpose of this note is to understand the structure and the orders of
the elements of these groups by means of the reduced $KO$-relations and the
reduced $J$-relations with as little elementary number theory as possible.
In order to achieve this, since the collection of all relations are very
complicated, we will slickly look at the effect of these relations on the
Atiyah-Hirzebruch spectral sequence (AHSS), [5] and do some counting which
helps us to understand the $J$-Group.

We will consider the problem also from a different perspective through the
rational representations of cyclic groups.

For simplicity, we will consider only lens spaces modulo a prime power $%
p^{r},$ $r\geq 1,$ where the prime $p$ is odd, throughout the paper. The
story is not so much different for the prime number $2$, [2].

\section{KO-Ring}

To avoid the complication of the top cell, let $L=L^{2k}$ be the $2k$-th
dimensional skeleton of the standard lens space of dimension $2k+1$ modulo $%
p^{r}$ and thus also of the classifying space $BZ_{p^{r}}$ of the cyclic
group $Z_{p^{r}}$. This space is our main object in this paper. Let $\eta $
denote the tautological complex line bundle over $BZ_{p^{r}}$ and over its
all skeletons and $\mu =\eta -1$ denote its reduction. Let $w=r(\mu )$ be
the realification of $\mu .$ Let $\psi ^{i}$ denote the real Adams operation
of degree $i\geq 1.$ Then, $\psi ^{i}$ acts on $w$ by the quadratic binomial
as $\psi ^{i}(w)=\dsum\limits_{j=1}^{i}\frac{\binom{i}{j}\binom{i+j-1}{j}}{%
\binom{2j-1}{j}}w^{j}.$ For odd $n,$ let $f_{n}(w)$ be the polynomial
defined as $f_{n}(w)=\psi ^{\frac{n+1}{2}}(w)-\psi ^{\frac{n-1}{2}}(w).$

$KO$-ring of $L$ is described in the following theorem.

\textbf{Theorem 1. }%
\begin{equation*}
KO(L)=Z[w]\diagup \left( f_{p^{r}}(w),w^{\left[ \frac{k}{2}\right] +1}\right)
\end{equation*}%
where $\left[ \cdot \right] $ is the floor function.

The relation $w^{\left[ \frac{k}{2}\right] +1}=0$ will be called the
terminating relation. The main relation $f_{p^{r}}(w)=0$ creates the
following set of relations by "complete reduction", [6], and by
multiplication by powers of $w$:

\textbf{Lemma 2. }$p^{r}w^{n}=-p^{r-1}w^{n+\frac{p-1}{2}}+h.o.t.$ for all $%
n\geq 1.$

In the above lemma, the higher order terms mean not only in the filtration
sense, but in the sense that multiplication of them by $p$ exterminates them
faster than $p^{r-1}w^{n+\frac{p-1}{2}}.$ Although, the set of relations in
Lemma 2 together with the terminating relation \ describe the additive
structure of the $KO$-ring, the decomposition of the reduced $\widetilde{KO}$%
-group into cyclic groups seems to be complicated and we skip this
discussion. But, ironically we will ask for the decomposition of the $%
\widetilde{J}$-group in the next section.

The main relation in the form of the set of the relations in Lemma 2 show
the effect of multiplication by $p$ on the filtrations of the AHSS of $KO(L)$
and this is the most essential observation for practical purposes. For
example, from the lemma, we can immediately calculate the order of the $%
\widetilde{KO}$-ring as $p^{r\left[ \frac{k}{2}\right] }$ which is also
clear from the AHSS and the $KO$-order of $w^{n}$ as $p^{r+\left[ \frac{k-2n%
}{p-1}\right] }$, [5].

\section{J-Group}

Let $q$ be a generator of the unit group $(Z\diagup p^{r})^{\star }.$ It
follows from the Adams conjecture that, [5],

\textbf{Theorem 3. }$J(L)=KO(L)\diagup (\psi ^{q}-1)KO(L).$

Note that the above quotient is the group quotient of the additive group $%
KO(L)$ by its subgroup $(\psi ^{q}-1)KO(L)$ which is obtained by the action
of the operation $\psi ^{q}-\psi ^{1}$. Hence, the set of relations that
defines $J(L)$ is $\psi ^{q}(w^{n})-w^{n}=0,$ $n\geq 1$. Note that these
relations are additive relations unlike the $KO$-relations.

Now let $I_{s}=\left\{ n\mid n=m\dfrac{p^{s}(p-1)}{2}\text{ where }\gcd
(m,p)=1\text{ and }n\leq \left[ \frac{k}{2}\right] \right\} $ for $s\geq 0$.
Let also $v_{p}(\cdot )$ denote the $p$-exponent function.

\textbf{Lemma 4. }In $J(L)$ we have

i) $p^{s+1}w^{n}=p^{s}w^{n+1}+h.o.t.$ for $n\in I_{s}$ and $0\leq s\leq r-1.$

ii) $w^{n}$ can be expressible by higher order terms if $n\notin I_{s}$ for
any $s.$

iii) $KO$-relations dominate for $n\in I_{s}$ and $s\geq r.$

\textit{Proof: }Since $\psi ^{q}(w^{n})=\left[ \psi ^{q}(w)\right] ^{n}$,
the relation $\psi ^{q}(w^{n})-w^{n}=0$ explicitly can be written in the
form 
\begin{equation*}
(q^{2n}-1)w^{n}=-\dfrac{q^{2n}(q^{2}-1)n}{12}w^{n+1}+h.o.t.
\end{equation*}%
Now, the lemma follows from the fact that $v_{p}(q^{2n}-1)=s+1$ if $n\in
I_{s}$ and $v_{p}(q^{2n}-1)=1$ otherwise and from Lemma 2.

Let $\widetilde{J}(L)$ denote the reduced $J$-group. It is a $p$-group being
a quotient of the $p$-group $\widetilde{KO}(L).$ From Lemma 4, we deduce the
following two interesting results which are also stated, for example, in
[7],[8] and [3], respectively.

\textbf{Theorem 5. }$v_{p}(\left\vert \widetilde{J}(L)\right\vert
)=\dsum\limits_{s=0}^{r-1}\left[ \dfrac{k}{p^{s}(p-1)}\right] .$

\textit{Proof:} Due to Lemma 2 (ii), if $n\notin I_{s}$ for any $s,$ then
the element $w^{n}$ and its multiples shouldn't be counted. Due to Lemma 2
(i), if $n\in I_{s}$ where $0\leq s\leq r-1,$ then the element $w^{n}$
contributes $p^{s+1}$ factor to $\left\vert \widetilde{J}(L)\right\vert $.
And due to Lemma 2 (iii), for other values of $n\leq \left[ \frac{k}{2}%
\right] ,$ the element $w^{n}$ contributes $p^{r}$ factor which is
independent of $s$. It follows from these observations that \textbf{\ }$%
v_{p}(\left\vert \widetilde{J}(L)\right\vert
)=\dsum\limits_{s=0}^{r-1}\left\vert I_{s}\right\vert $. But, it is not
difficult to see that $\left\vert I_{s}\right\vert =\left[ \frac{k}{%
p^{s}(p-1)}\right] $.

\textbf{Theorem 6. }Let $\left\vert w\right\vert $ denote the $J$-order of
the element $w$. Then 
\begin{equation*}
v_{p}(\left\vert w\right\vert )=\max \left\{ s+\left[ \frac{k}{p^{s}(p-1)}%
\right] p^{s}\mid 0\leq s\leq r-1\text{ and }p^{s}(p-1)\leq k\right\} .
\end{equation*}

\textit{Proof:} If $n$ is the greatest number for which $n\in I_{s}$ for a
fixed $0\leq s\leq r-1,$ then we have $\left\vert w\right\vert \geq s+\frac{%
2n}{p-1}=s+\left[ \frac{k}{p^{s}(p-1)}\right] p^{s}.$ Therefore, we should
take the maximum over all $s$'s.

Now, we turn to the main question: The group structure of $J(L).$ In other
words, we want to determine the (natural) decomposition of the finite
Abelian group $\widetilde{J}(L)$ into cyclic groups. Let $Z_{a}(g),$ $Z(g)$
denote the cyclic group $G$ with the generator $g$ of order $a$ and of order
infinity respectively$.$

\textbf{Theorem 7.}$\mathbf{\ }\widetilde{J}(BZ_{p^{r}})=\dbigoplus%
\limits_{t=0}^{r-1}Z(\psi ^{p^{t}}(w)).$

\textit{Proof:} We will prove by induction on $r.$ For $r=1,$ $\psi
^{p}(w)=0.$ Hence, it follows from Theorem 3 that multiples of $w$ span $%
\mathbf{\ }\widetilde{J}(BZ_{p^{r}}).$ Now, let it be true for $r-1$, that
is for the space $BZ_{p^{r-1}}$. The $p$-th power homomorphism, that is the
transfer $\tau :\mathbf{\ }\widetilde{J}(BZ_{p^{r}})\rightarrow \mathbf{\ }%
\widetilde{J}(BZ_{p^{r-1}})$ is onto with kernel $Z(\psi ^{p^{r-1}}(w))$.

The group $\widetilde{J}(L)$ is clearly a quotient group of $\widetilde{J}%
(BZ_{p^{r}})$ because of Theorem 3. Unfortunately things are very
complicated in this quotient, [1], [8]. On the other hand, we can find an
upper bound for the group $\widetilde{J}(L)$ by determining the order of the
elements $\psi ^{p^{s}}(w)$. The counting here is generalization of Theorem
6.

\textbf{Theorem 8. }$v_{p}(\left\vert \psi ^{p^{t}}(w)\right\vert )=\max
\left\{ s-t+\left[ \dfrac{k}{p^{s}(p-1)}\right] p^{s-t}\mid t\leq s\leq r-1%
\text{ and }p^{s}(p-1)\leq k\right\} .$

You can see the above order in [3] as well. The following is a corollary of
Theorems 7 and 8.

\textbf{Theorem 9. }Let $N=\min \left\{ r-1,\left[ \log _{p}(k+1)\right]
\right\} $ and $a_{t}=\left\vert \psi ^{p^{t}}(w)\right\vert .$ Then 
\begin{equation*}
\widetilde{J}(L)=\dbigoplus\limits_{t=0}^{N}Z_{a_{t}}(\psi
^{p^{t}}(w))\diagup I
\end{equation*}%
where $I$ is the subgroup generated by all linear combinations of $\left\{
w,\psi ^{p}(w),\psi ^{p^{2}}(w),...,\psi ^{p^{N}}(w)\right\} $ which are
zero in $\widetilde{J}(L).$

The bothering group $I$ occurs because of the terminating relation $w^{\left[
\frac{k}{2}\right] +1}=0$ as it happened in the calculation of the number $%
a_{t}$. Relatively, $I$ is a very small subgroup compared to $%
\dbigoplus\limits_{t=0}^{N}Z_{a_{t}}(\psi ^{p^{t}}(w))$ itself. But, it is
responsible for all complications and the reason behind the messy work [1].
This is also the case for the works [8] and [2]. It is even not clear
whether the results in these articles match. But, in some special cases,
this group is zero. For example,

\textbf{Corollary 10.} Let $r=1.$ Then $\widetilde{J}(L)=Z_{a_{0}}(w)$ where 
$a_{0}=p^{[\frac{k}{p-1}]}.$

In the next section, we will consider the same problem from a different
perspective.

\section{Rational Representations}

Topological $J$-Theory has a mirror image in representation theory like
topological $K$-theory through Atiyah-Segal completion theorem isomorphism,
[4]. It is interesting that $J(BZ_{p^{r}})$ has some deep connections with
the rational representations of the cyclic group $Z_{p^{r}}.$

Let $L$ and $q$ be as described in the previous sections. Let $KQ(L)$ denote
the kernel of the group homomorphism $K(L)\overset{\psi ^{q}-1}{\rightarrow }%
K(L).$ Then, we have an (additive) group isomorphism $\widetilde{J}(L)\cong 
\widetilde{KQ}(L)$ since kernel and cokernel are isomorphic for an ordered
group endomorphism of a finite Abelian group -it is not for any
endomorphism- and the complex and the real $J$-groups of $L$ are the same.
We will call $KQ(L)$ as the rational topological $K$-theory of $L$. The
reason of that name is justified in Theorem 12 below.

First, we will define some essential reduced elements in $KQ(L).$ Recall
that $\eta $ is the Hopf bundle.

Definition 11\textbf{. }Let $\Omega _{0}=\left( \dsum\limits_{s^{\ast }=1%
\text{ }}^{p^{r}-1}\eta ^{s}\right) -(p^{r}-p^{r-1})$ where $\gcd (s,p)=1$
in the summation. And let $\Omega _{t}$ be defined inductively as $\Omega
_{t}=\frac{1}{p}\psi ^{p}(\Omega _{t-1})$ for $1\leq t\leq r-1.$

Note that these are reduced vector bundles corresponding to the rational
representations of the group $Z_{p^{r}}$. The following is the analogue of
Theorem 9.

\textbf{Theorem 12. }$\widetilde{KQ}(L)=\dbigoplus%
\limits_{t=0}^{N}Z_{a_{t}}(\Omega _{t})\diagup I$ where $N$ and $a_{t}$ are
as defined in Theorem 9 and 8 and where $I$ is the subgroup generated by all
linear combinations of $\left\{ \Omega _{0},\Omega _{1},...,\Omega
_{N}\right\} $ which are zero in $\widetilde{KQ}(L)$.

\textit{Proof:} It can be proved by induction on $r$ that $\widetilde{KQ}%
(BZ_{p^{r}})$ is direct sum of infinite cyclic subgroups whose generators
are $\Omega _{t}$'s $1\leq t\leq r-1$ as in Theorem 7. $\widetilde{KQ}(L)$
is a quotient group of this and the orders of $\Omega _{t}$ in this quotient
must match with Theorem 9 according to the group isomorphism mentioned above.

It is interesting that $KQ(L)$ seems to have a ring structure unlike the
cokernel $J(L)$. We can derive the analogue of Corollary 10 and actually
more of it independently by using the multiplicative structure of the ring $%
K(L)$.

\textbf{Theorem 13. }Let $r=1.$ Then $\widetilde{KQ}(L)=Z[\Omega
_{0}]/(\Omega _{0}^{2}+p\Omega _{0},\Omega _{0}^{\left[ \frac{k}{p-1}\right]
+1})$.

\textit{Proof:} For $r=1,$ $\Omega _{0}=\left( \dsum\limits_{s=1\text{ }%
}^{p-1}\eta ^{s}\right) -(p-1).$ Now, if we square this in $K(L)$ and use $%
\eta ^{p}=1$ for simplifications, we get $\Omega _{0}^{2}=-p\Omega _{0}.$The
filtration of $\Omega _{0}$ on the AHSS of $K(L)$ is same as $\mu ^{p-1}.$
Since $\Omega _{0}=\frac{(1+\mu )^{p}-1-p\mu }{\mu }$, the filtration of $%
\Omega _{0}^{s}$ is same as $\mu ^{s(p-1)}$. Now, $\mu ^{s(p-1)}=0$ in $K(L)$
when $s(p-1)>k$, [5]. Hence, when $s>\left[ \frac{k}{p-1}\right] $, $\Omega
_{0}^{s}=0$ in $\widetilde{KQ}(L)$. Thus, we also get the terminating
relation $\Omega _{0}^{\left[ \frac{k}{p-1}\right] +1}=0.$

It easily follows from the above theorem that $\widetilde{KQ}(L)$ is a
cyclic group with the generator $\Omega _{0}$ and order $p^{[\frac{k}{p-1}]}$
which agrees with Corollary 10. It is also interesting to compare \ the
complicated $J$-relations $\psi ^{q}(w^{n})-w^{n}=0$ in previous approach
with the simple multiplicative relation $\Omega _{0}^{2}+p\Omega _{0}=0$ in
Theorem 13. It seems that these kind of relations are valid for the general
case. In $KQ(L),$ it seems that we have $\Omega _{t}^{2}=-p\Omega _{t}$ for
any $1\leq t\leq r-1.$ It remains only to understand how the terminating
relation $\mu ^{k+1}=0$ turns to a relation in terms of $\Omega _{t}$'s,
illuminating the problematic tiny subgroup $I.$ Thus, we have come to the
same dead-lock of the previous section, but now the relations here are
simple and multiplicative. They are more promising.

\end{document}